\documentclass[a4paper,12pt]{article}

\usepackage{amsmath,amssymb,amsthm}
\usepackage{enumerate}
\usepackage{epic,eepic}

\newtheorem{thm}{Theorem}
\theoremstyle{remark}
\newtheorem*{ohno}{Ohno relation}
\theoremstyle{plain}
\newtheorem{lem}{Lemma}
\newtheorem{prop}{Proposition}
\newtheorem{cor}{Corollary}
\theoremstyle{definition}
\newtheorem{defn}{Definition}

\newcommand{\Z}{\mathbb Z}
\newcommand{\Q}{\mathbb Q}

\newcommand{\bfc}{{\boldsymbol c}}
\newcommand{\bfk}{{\boldsymbol k}}

\newcommand{\Li}[2]{\mathrm{Li}_{#1}\left(#2\right)}
\newcommand{\imu}{\sqrt{-1}}

\DeclareMathOperator*{\summ}{{\sum}^\prime}
\DeclareMathOperator*{\Res}{\mathrm{Res}}
\DeclareMathOperator{\len}{\mathrm{len}}

\allowdisplaybreaks[2]

\title{
  Relations for Multiple Zeta Values\\
  and\\
  Mellin Transforms of
  Multiple Polylogarithms
}
\author{Jun-ichi OKUDA \footnote{E-mail:  \tt{okuda@gm.math.waseda.ac.jp}}
 \qquad Kimio UENO \footnote{E-mail:  \tt{uenoki@mse.waseda.ac.jp}}\\
  \\
  Department of Mathematical Sciences,\\
  School of Science and Engineering,\\
  WASEDA UNIVERSITY,\\
  3-4-1, Okubo Shinjuku-ku,
  Tokyo 169-8555, Japan
}

\date{}

\begin{document}
\maketitle

\begin{abstract}
 In this paper a relationship
 between the Ohno relation for multiple zeta values
 and multiple polylogarithms are discussed.
 First we introduce generating functions for the Ohno relation, and
 investigate their properties. We show that there exists a subfamily
 of the Ohno relation which recovers algebraically its totality.
 This is proved through analysis of Mellin transform of
 multiple polylogarithms. Furthermore, this subfamily is shown
 to be converted to the Landen connection formula for multiple
 polylogarithms by inverse Mellin transform.
\end{abstract}

\section{Introduction}

\subsection{Definitions and examples}

In this paper, we will consider the relationship
between the Ohno relation for multiple zeta values
(MZVs, for short)
and the Landen connection formula for multiple polylogarithms
(MPLs, for short)
via Mellin transform and inverse Mellin transform.

\begin{defn}
For positive integers $k_1,\dots,k_n$ and $|z|<1$,
MPLs are defined by
\begin{align}
 \Li{k_1,k_2,\dots,k_n}{z}
 &:= \sum_{m_1>m_2>\dots>m_n>0}
    \frac{z^{m_1}}{m_1^{k_1}m_2^{k_2}\dotsm m_n^{k_n}}
\end{align}
and for the null sequence $\emptyset$, $\Li{\emptyset}{z} := 1$.
If $k_1\ge2$,  MPLs also converge at $z=1$ and define MZVs
\begin{align}
\zeta(k_1,k_2,\dots,k_n)
 &:=
 \sum_{m_1 > m_2 > \dots >m_n>0}
  \frac{1}{m_1^{k_1}m_2^{k_2}\dotsb m_n^{k_n}},
\end{align}
and similarly $\zeta(\emptyset) := 1$.
The weight and the depth of $\zeta(k_1,\dots,k_n)$ are 
defined to be 
$k_1 + \dotsb + k_n$ and $n$,  respectively.
\end{defn}

Through consideration on the dilogarithm $\Li{2}{z}$ 
as an example, we explain an essential aspect of 
the relationship that will be considered in our papers.

The sum formula \cite{gran}
for MZVs of depth $2$
\begin{align*}
  \zeta(3+l)=
  \sum_{\substack{c_1+c_2=l\\c_1,c_2\ge0}}\zeta(2+c_1,1+c_2)
  \qquad (l\in\Z_{\ge0})
\end{align*}
is equivalent to the generating functional expression
\begin{align*}
 \sum_{l=0}^\infty \zeta(3+l)\,\lambda^l
  = \sum_{l=0}^\infty
    \left\{
     \sum_{\substack{c_1+c_2=l\\c_1,c_2\ge0}}\zeta(2+c_1,1+c_2)
    \right\}
     \lambda^l.
\end{align*}
Noting that, for a positive integer $n$, 
\begin{align*}
\frac{1}{n-\lambda} = \sum_{l=0}^{\infty}\frac{\lambda^l}
                  {n^{l+1}} \qquad \text{for} \quad |\lambda|<1,
\end{align*}
we see that
the both sides in the above are meromorphic functions in $\lambda$,
\begin{align*}
  \sum_{n=1}^\infty \frac{1}{n^2(n-\lambda)} 
  =
 \sum_{n_1>n_2>0} \frac{1}{n_1(n_1-\lambda)(n_2-\lambda)}.
\end{align*}
Applying ``inverse Mellin transform''
\begin{align*}
  \widetilde{M}\left[f(\lambda)\right](z)
   = \frac{1}{2\pi \imu}
       \int_C
        f(\lambda)z^\lambda\ d\lambda
  \qquad(0<z<1)
\end{align*}
to the left hand side
(for details see Section~\ref{sec.mellin}),
we have
\begin{align*}
 \sum_{n=1}^\infty \frac{z^n}{n^2}
 =\Li{2}{z}.
\end{align*}
Note that
\begin{multline*}
  \sum_{n_1>n_2>0} \frac{1}{n_1(n_1-\lambda)(n_2-\lambda)}\\
  =
  \sum_{n_1>n_2>0} \frac{1}{n_1(n_1-n_2)(n_2-\lambda)}
 +\sum_{n_1>n_2>0} \frac{1}{n_1(n_2-n_1)(n_1-\lambda)}.
\end{multline*}
Applying inverse Mellin transform to each term,
we have
\begin{align*}
 \sum_{n_1>n_2>0} \frac{z^{n_2}}{n_1(n_1-n_2)}
 &= -\Li{2}{\frac{z}{z-1}},
 \\
 \sum_{n_1>n_2>0} \frac{z^{n_1}}{n_1(n_2-n_1)}
 &= - \Li{11}{\frac{z}{z-1}}.
\end{align*}
Consequently we obtain the next functional equation
for the dilogarithm:
\begin{align}
 \Li{2}{z} = -\Li{2}{\frac{z}{z-1}} - \Li{11}{\frac{z}{z-1}},
 \label{eq.landilog}
\end{align}
which is known as the Landen connection formula
for the dilogarithm \cite{lewin}.
This can be viewed as the connection formula for the dilogarithm
between $1$ and $\infty$.

\subsection{Main results and organization}
Now we explain the Ohno relation \cite{ohno}
which is a generalization of the sum formula.

Any index $\bfk = (k_1,\dots,k_n) \in \Z_{\ge 1}^n$, $k_1\ge 2$
can be written uniquely as 
\begin{align}
 \bfk &= (a_1+1,\underbrace{1,\dots,1}_{b_1-1},\dots,
  			a_s+1,\underbrace{1,\dots,1}_{b_s-1}),\\
 \intertext{with $s\in\Z_{\ge 1}$
  and $a_i,b_i\in\Z_{\ge 1}\,(i=1,\dotsc,s)$.
  The dual index $\bfk' = (k_1',\dots,k_{n'}')$ of $\bfk$ is defined by}
  \bfk'  &= (b_s+1,\underbrace{1,\dots,1}_{a_s-1},\dots,
  			b_1+1,\underbrace{1,\dots,1}_{a_1-1}),
  \label{eq.dual}
\end{align}
and the dual of $\emptyset$ is itself.

\begin{ohno}\label{thm.ohno}
Let $\bfk=(k_1,\dots,k_n)$ be any index and $\bfk'$ be its dual.
For all $l\in\Z_{\ge 0}$,
we have the following homogeneous (w.r.t.\ weight) relation,
\begin{align}
  \sum_{\substack{c_1 + \dots + c_n = l\\c_j\ge0}}
   \!\!
   \zeta(k_1+c_1,\dots,k_n+c_n)
  =
   \!\!
  \sum_{\substack{c_1' + \dots + c_{n'}' = l\\c_j'\ge0}}
   \!\!
   \zeta(k_1'+c_1',\dots,k_{n'}'+c_{n'}').
 \end{align}
In particular,
this contains such relations
as the Hoffman relation \cite{hoff1} ($l=1$),
the duality formula \cite{zagi} ($l=0$)
and the sum formula ($n=1$).
\end{ohno}

We introduce the generating functions
of the both sides of the Ohno relation
as follows:
\begin{align}
 f((a_i,b_i)_{i=1}^s;\lambda)
 &:=
 \sum_{l=0}^{\infty}
 \left\{
  \sum_{\substack{c_1+\dots+c_n=l\\c_j\ge0}}
   \zeta(k_1+c_1,\dots,k_n+c_n) \right\}
 \lambda^l,\\
 g((a_i,b_i)_{i=1}^s;\lambda)
 &:=
 \sum_{l=0}^{\infty}
 \left\{
  \sum_{\substack{c'_1+\dots+c_{n'}'=l\\c'_j\ge0}}
   \zeta(k'_1+c'_1,\dots,k_{n'}'+c_{n'}') \right\}
 \lambda^l.
\end{align}
The Ohno relation reads
\begin{align}
 f((a_i,b_i)_{i=1}^s;\lambda)
 =
 g((a_i,b_i)_{i=1}^s;\lambda).
 \label{eq.func ohno}
\end{align}
We can show that $f$'s and $g$'s satisfy
the same difference relations
which play a fundamental role in our theory
(Proposition~\ref{prop.sabun}, Section~\ref{sec.diffrel}).
Now we define the functions
$F(\bfk;\lambda)$ and $G(\bfk;\lambda)$ by
\begin{align*}
  F(\bfk;\lambda) &:=
   \sum_{\delta_i=0,1} (-\lambda)^{n-1-|\delta|}
   f((k_1,1)\cup(k_i-\delta_i,1)_{i=2}^n;\lambda),\\
  G(\bfk;\lambda) &:=
   \sum_{\delta_i=0,1} (-\lambda)^{n-1-|\delta|}
   g((k_1,1)\cup(k_i-\delta_i,1)_{i=2}^n;\lambda).
\end{align*}
Under the Ohno relation, we have
\begin{align}
 F(\bfk;\lambda)=G(\bfk;\lambda),
 \label{eq.reduced ohno}
\end{align}
which we call the reduced Ohno relation.
We can easily see that
\begin{align*}
   F(\bfk;\lambda)=\sum_{m_1>m_2>\dots>m_n>0}
      \frac{1}{m_1^{k_1}(m_1-\lambda)m_2^{k_2}
                       \cdots m_n^{k_n}},
\end{align*}
so that inverse Mellin transform of $F(\bfk;\lambda)$ is
the multiple polylogarithm $\Li{\bfk}{z}$.
This fact gives us strong motivation to 
introduce these functions. The main theorem of this paper
(Theorem~\ref{thm.redohno}) 
says that
{\itshape the generating functions $f$'s and $g$'s
 are represented as
MZVs-linear combinations of $F$'s and $G$'s, respectively.
In other words, the reduced Ohno relation \eqref{eq.reduced ohno} recovers
the totality of the Ohno relation \eqref{eq.func ohno}}
(Section~\ref{sec.redohno}).
This theorem is proved by virtue of the differential equations
satisfied by MPLs and Mellin transform
(Section~\ref{sec.mellin}~and~\ref{sec.proofredohno}).
Furthermore via inverse Mellin transform
{\itshape the reduced Ohno relation is converted to
the Landen connection formula for MPLs$:$
\begin{align}
   \Li{k_1,\dots,k_n}{z} =
   (-1)^n \sum_{\substack{c_1,\dots,c_n\\|c_j|=k_j}}
    \Li{c_1\dotsm c_n}{\frac{z}{z-1}},
    \label{eq.landen}
\end{align}
where $c_i$ runs all compositions of $k_i$
and the product of $c_j$'s is given by concatenation}.
This can be viewed as the connection formula 
between $1$ and $\infty$ for MPLs
(Section~\ref{sec.landen}).
In Section~\ref{sec.proofsabun}
we give the proof of Proposition~\ref{prop.sabun}
and another proof of the Ohno relation.

\section*{Acknowledgements}

 The authors express their deep gratitude
 to Professor Masanobu Kaneko and Professor Shigeki Akiyama
 for their valuable suggestions.

 The second author is partially supported by
 Grant-in-Aid Scientific Research from the Ministry of 
 Education, Culture, Sports, Science and Technology 
 of Japan (12640046)
 and by Waseda University Grant for Special Research Project
 (2000A-124, 2001A-088 and 2002A-067 ).

\section{The generating functions and their properties}
\label{sec.diffrel}

\subsection{Compositions}

By a composition of a positive integer $n$, we mean
an ordered sequence $\bfc = (c_1,\dots,c_l)$
of positive integers of which the sum is equal to $n$,
and the composition of $0$ is defined to be $\emptyset$.
The ``weight'' $|\bfc|$ and the ``length'' $\len(\bfc)$ of $\bfc$ are,
by definition, 
$n$ and $l$ respectively.
We allow $0$'s to appear in the middle elements of $\bfc$
and identify such compositions and normal compositions
by removing $0$'s successively,
i.e.\ we regard $(\dots,\underline{c_{i-1},0,c_{i+1},}\dots)$
to be the same composition as
$(\dots,c_{i-1}+c_{i+1},\dots)$.
For example,
\begin{align*}
 (\underline{3,0,2},0,4) &= (\underline{5,0,4}) = (9),
 \quad\text{or}\quad
 (3,0,\underline{2,0,4}) = (\underline{3,0,6}) = (9),\\
 (3,0,2,1,\underline{0,0,4})
  &= (\underline{3,0,2},1,4)= (5,1,4).
\end{align*}
We remark that
the resulting composition dose not depend on
the procedure of the identification.
For compositions $\bfc$ and $\bfc'$
we define the partial order $\bfc\succ \bfc'$
if $\bfc'$ is obtained by
decreasing some elements of $\bfc$.
For example, $(5,1,4) \succ (4,1,3)$
and $(5,1,4) \succ (9) = (5,0,4)$.
Between compositions of even length
and compositions whose first element is greater than $1$,
we define the 1-1 correspondence $\kappa$ as follows:
\begin{align}
 \kappa((a_i,b_i)_{i=1}^s)
 = (a_1+1,\underbrace{1,\dots,1}_{b_1-1},
    \dots,a_s+1,\underbrace{1,\dots,1}_{b_s-1}).
\end{align}

\subsection{Generating functions for the Ohno relation}

\begin{defn}\label{defn.genfun}
For any composition $(a_i,b_i)_{i=1}^s=\kappa^{-1}(k_1,\dots,k_n)$,
we set the generating functions of MZVs as
\begin{align}
 f((a_i,b_i)_{i=1}^s;\lambda)
 &:=
 \sum_{l=0}^{\infty}
 \left\{
  \sum_{\substack{c_1+\dots+c_n=l\\c_j\ge0}}
   \zeta(k_1+c_1,\dots,k_n+c_n) \right\}
 \lambda^l,\\
 g((a_i,b_i)_{i=1}^s;\lambda)
 &:=
 \sum_{l=0}^{\infty}
 \left\{
  \sum_{\substack{c'_1+\dots+c_{n'}'=l\\c'_j\ge0}}
   \zeta(k'_1+c'_1,\dots,k_{n'}'+c_{n'}') \right\}
 \lambda^l\\
 &= f((b_i,a_i)_{i=s}^1;\lambda),\notag
\end{align}
where $(k_1',\dots,k_{n'}')$ is the dual of $(k_1,\dots,k_n)$.
For convenience  $f((a_i,b_i)_{i=1}^s;\lambda) := 0$
if $a_1$ or $b_s = 0$.
We set the weight of $f((a_i,b_i)_{i=1}^s;\lambda)$
and $g((a_i,b_i)_{i=1}^s;\lambda)$ to be
$|(a_i,b_i)_{i=1}^s|$.
\end{defn}

This power series absolutely converges for $|\lambda| < 1$.
The Ohno relation reads
\begin{align}
 f((a_i,b_i)_{i=1}^s;\lambda) = g((a_i,b_i)_{i=1}^s;\lambda),
\end{align}
for any compositions $(a_i,b_i)_{i=1}^s$.

Noting that for a positive integer $n$, 
\begin{align*}
\frac{1}{n-\lambda} = \sum_{l=0}^{\infty}\frac{\lambda^l}
                  {n^{l+1}} \qquad \text{for} \quad |\lambda|<1,
\end{align*}
one can easily see that 
\begin{align}
 f((a_i,b_i)_{i=1}^s;\lambda)
 =
 \sum_{m_1 > \dots > m_{B_s > 0}}
   \prod_{i=1}^s
    \frac{1}{m_{B_{i-1}+1}^{a_i}
     \underbrace{(m_{B_{i-1}+1}-\lambda)\dotsm(m_{B_i}-\lambda)}_{b_i}
    },
\end{align}
where $B_0=0$ and $B_i = b_1+\dotsb+b_i$ for $i\ge1$.

\subsection{Properties of the generating functions}

The generating functions satisfy the following 
difference equations:

\begin{prop}\label{prop.sabun}
 We set
 $\lambda' := \lambda - 1$ and 
 $I:=\left\{(0,0),(1,0),(0,1)\right\}$,
 then for any composition $(a_i,b_i)_{i=1}^s$
 the generating function $f$ satisfies the following relations.
  \begin{multline}
    \sum
    (-\lambda)^{s-|\delta|-|\epsilon|}
    f((a_i-\delta_i,b_i-\epsilon_i)_{i=1}^s;\lambda)\\
   = 
   \summ
    (-\lambda')^{s-|\delta'|-|\epsilon'|}
     f((a_i-\delta_i',b_i-\epsilon_{i+1}')_{i=1}^s;\lambda').
  \end{multline}
 Here the sum $\sum$ is taken over $\{\delta_i,\epsilon_i\}\in I$,
 the sum $\summ$ taken over $\delta_1', \epsilon_{m+1}'\in \{0,1\}$
 and $\{\delta_i', \epsilon_i'\}\in I$ for $i=2,\dots,m$,
 and $|\delta^{(\prime)}|$ (resp.\ $|\epsilon^{(\prime)}|$) is
 the sum of all $\delta_i^{(\prime)}$ (resp.\ $\epsilon_i^{(\prime)}$).
 The generating function $g$ also satisfies the same relations.
 We define the weight of $\lambda$ and $\lambda'$ to be $-1$
 and this relation is homogeneous of weight $|(a_i,b_i)_{i=1}^s| - s$.
\end{prop}

The proposition will play a crucial role in our theory.
The proof is so long
that it will be postponed until Section~\ref{sec.proofsabun}.

The generating functions are analytically continued to
meromorphic functions with simple poles at positive integers.

\begin{prop}\label{prop.partial}
The generating function $f((a_i,b_i)_{i=1}^s;\lambda)$
can be expanded to a partial fraction
 \begin{align}
  f((a_i,b_i)_{i=1}^s;\lambda)=
  \sum_{p=1}^\infty
   \left\{
    \sum_{j=1}^{B_s}
    \sum_{\substack{m_1 > \dots > m_{j-1} > p\\ p >m_{j+1}>\dots>m_{B_s}}}
     C_p^{m_1\dotsc \stackrel{\stackrel{j}\smallsmile}p \dotsc m_{B_s}}\right\}
    \frac{1}{p-\lambda}
\end{align}
where
\begin{align}
C_{m_j}^{m_1\dotsc m_{B_s}}
  = \frac{1}{m_1^{a_1}m_{B_1+1}^{a_2}\dotsc m_{B_{s-1}+1}^{a_s}}
          \prod_{i\not=j}\frac{1}{(m_i-m_j)}.
 \end{align}
\end{prop}

\begin{proof}
The generating function $f$ can be written as follows:
\begin{align*}
f((a_i,b_i)_{i=1}^s;\lambda)
 = 
  \sum_{m_1 > \dots > m_{B_s} > 0}
  \sum_{j=1}^{B_s}
   \frac{C_{m_j}^{m_1\dotsc m_{B_s}}}{m_j-\lambda}.
\end{align*}
For the proof, we have to show that it is possible to
change the order of the summations.
So it is sufficient to prove that for any $j$
\begin{align*}
 \sum_{\substack{m_1 > \dots > m_{j-1} > m_j\\ m_j >m_{j+1}>\dots>m_{B_s}}}
   \frac{C_{m_j}^{m_1\dotsc m_{B_s}}}{m_j-\lambda}
\end{align*}
converges absolutely.
Put $d_i = m_i - m_{i+1}$ for $i=1,\dotsc,B_s-1$ and $d_{B_s}= m_{B_s}$.
Making use of the inequality
\begin{align*}
 d_1 + d_2 + \dotsb + d_{B_s} \ge B_s \root{B_s}\of{d_1d_2\dotsb d_{B_s}}
\end{align*}
we have
\begin{align*}
 \left|\frac{C_{m_j}^{m_1\dotsc m_{B_s}}}{m_j-\lambda}\right|
 &=\left|
   \left(
     \frac{1}{m_1^{a_1}\dotsm m_{B_{s-1}+1}^{a_s}}
          \prod_{i\not=j}\frac{1}{(m_i-m_j)}\right)
    \frac{1}{m_j-\lambda}
  \right|\\
 &\le
  \frac{1}{(d_1+d_2+\dotsb+d_{B_s})^{a_1}}
   \left(\prod_{i=1}^{B_s-1}\frac{1}{d_i}\right)
    \frac{1}{\left|d_j+\dotsb+d_{B_s}-\lambda\right|}\\
 &\le
  \frac{1}{\left(B_s \root{B_s}\of{d_1d_2\dotsb d_{B_s}}\right)^{a_1}}
   \left(\prod_{i=1}^{B_s-1}\frac{1}{d_i}\right)
    \frac{1}{\left|d_j+\dotsb+d_{B_s}-\lambda\right|}
\end{align*}
Let $\lambda$ be in a compact set
which dose not involve positive integers.
Then there exists a positive constant $A$ such that
\begin{align*}
 \frac{1}{\left|d_j+\dotsb+d_{B_s}-\lambda\right|}
 \le
 \frac{A}{d_{B_s}}.
\end{align*}
Hence
\begin{align*}
 \sum_{m_1 > \dots > m_{B_s} > 0}
  \left|\frac{C_{m_j}^{m_1\dotsc m_{B_s}}}{m_j-\lambda}\right|
 &\le
 A
 \sum_{d_1,\dotsc,d_{B_s}=1}^\infty
  \frac{1}{\left(B_s \root{B_s}\of{d_1\dots d_{B_s}}\right)^{a_1}}
   \left(\prod_{i=1}^{B_s}\frac{1}{d_i}\right)\\
 &<+\infty.
\end{align*}
\end{proof}

\section{Algebraic reduction of the Ohno relation}
\label{sec.redohno}

\begin{defn}\label{defn.FG}
 For any index $\bfk=(k_1,\dots,k_n)$
 we set the homogeneous functions of weight $|\bfk|+1$
 as
 \begin{align}
  F(\bfk;\lambda) &:=
   \sum_{\delta_i=0,1} (-\lambda)^{n-1-|\delta|}
   f((k_1,1)\cup(k_i-\delta_i,1)_{i=2}^n;\lambda),
   \label{eq.F}\\
  G(\bfk;\lambda) &:=
   \sum_{\delta_i=0,1} (-\lambda)^{n-1-|\delta|}
   g((k_1,1)\cup(k_i-\delta_i,1)_{i=2}^n;\lambda).
 \end{align}
\end{defn}

It is easy to calculate $F(\bfk;\lambda)$; we have
\begin{align}
   F(\bfk;\lambda)=\sum_{m_1>m_2>\dots>m_n>0}
      \frac{1}{m_1^{k_1}(m_1-\lambda)m_2^{k_2}
                       \cdots m_n^{k_n}}.
  \label{eq.seriesf}
\end{align}
On the other hand, it is difficult to write
down the explicit form of $G(\bfk;\lambda)$.
These functions satisfy difference equations of the simple
form:

\begin{prop}\label{prop.Frel}
 For any index $\bfk=(k_1,\dots,k_n)$, we have the relations 
 homogeneous of weight $|\bfk|$$:$
\begin{enumerate}[{\normalfont \upshape (i)}]
 \item if $k_1\ge2$\\
 $\bfk'$ to be the dual index of $\bfk$
 defined by \eqref{eq.dual}.
 Then
 \begin{align}
  \lambda F(k_1,k_2,\dots,k_n;\lambda) + \zeta(k_1,\dots,k_n)
  &=  F(k_1-1,k_2,\dots,k_n;\lambda),\\
  \label{eq.const}
  \lambda G(k_1,k_2,\dots,k_n;\lambda) + \zeta(k'_1,\dots,k'_{n'})
  &=  G(k_1-1,k_2,\dots,k_n;\lambda).
 \end{align}
 \item if $k_1=1$\\
 $(k'_2+1,k_3,\dots,k'_{n'})$ to be
 the dual index of $(k_2+1,k_3,\dots,k_n)$.
 Then
 \begin{align}
 \begin{split}
  &\lambda F(1,k_2,\dots,k_n;\lambda)
   + \zeta(k_2+1,k_3,\dots,k_n)\\
   &\qquad\qquad\qquad=
    \lambda' F(1,k_2,\dots,k_n;\lambda')
     +\lambda' F(k_2+1,\dots,k_n;\lambda'),
 \end{split}\\
 \begin{split}
  &\lambda G(1,k_2,\dots,k_n;\lambda)
   + \zeta(k'_2+1,k'_3,\dots,k'_{n'})\\
   &\qquad\qquad\qquad=
    \lambda' G(1,k_2,\dots,k_n;\lambda')
     +\lambda' G(k_2+1,\dots,k_n;\lambda').
 \end{split}
 \end{align}
\end{enumerate}
\end{prop}

\begin{proof}
 Induction on compositions.
 Apply Proposition~\ref{prop.sabun} to $(k_i,1)_{i=1}^n$
 and gather $f$'s or $g$'s whether $\epsilon_i = 0$ or not,
 then there are many cancel outs
 because of the identification
 $(\dots,k_i-1,0,k_{i+1},\dots) = (\dots,k_i,0,k_{i+1}-1,\dots)$
 and the induction hypothesis.
\end{proof}

It is easy to see that the inverse Mellin transform
of $F(\bfk;\lambda)$ is the MPL $Li_{\bfk}(z)$. This is a
motivation to introduce these functions.
It is known that the Ohno relation
is the largest systematic relation for MZVs,
however there are many linear dependency among them.
Acutually we can prove

\begin{thm}\label{thm.redohno}
For any composition $(a_i,b_i)_{i=1}^s$, we have
\begin{align}
  f((a_i,b_i)_{i=1}^s;\lambda)
  =
  \sum_{\bfc}
   \alpha_\bfc^{(a_i,b_i)} \zeta(\bfk_\bfc^{(a_i,b_i)})
    F(\bfc;\lambda),
  \label{eq.fasF}
  \\
  g((a_i,b_i)_{i=1}^s;\lambda)
  =
  \sum_{\bfc}
   \alpha_\bfc^{(a_i,b_i)} \zeta(\bfk'^{(a_i,b_i)}_\bfc)
    G(\bfc;\lambda),
\end{align}
where the summation runs over
some finite number of compositions $\bfc$,
$\alpha_{\bfc}^{(a_i,b_i)}\in\Q$,
$\kappa^{-1}(\bfk_{\bfc}^{(a_i,b_i)})) \prec (a_i,b_i)_{i=1}^s$,
and
$\bfk'^{(a_i,b_i)}_{\bfc}$ is the dual index
for $\bfk_{\bfc}^{(a_i,b_i)}$.
Moreover the duality formula
$\zeta(\bfk_{\bfc}^{(a_i,b_i)})=\zeta(\bfk'^{(a_i,b_i)}_{\bfc})$
comes from $f(0)=g(0)$
for compositions less than $(a_i,b_i)_{i=1}^s$.
So the Ohno relation is reduced to
\begin{align}
  F(\bfk;\lambda) = G(\bfk;\lambda)
  \label{eq.redohno}
\end{align}
as an algebraic relation.
\end{thm}

We call $(\ref{eq.redohno})$ the reduced 
Ohno relation. For the proof of the theorem,
we have to consider inverse Mellin trasform of $f$'s and $g$'s.

\section{Inverse Mellin transform of the generating functions}
\label{sec.mellin}

\subsection{Integral transform of the generating functions}
For any composition $(a_i,b_i)_{i=1}^s$
and any integer $l$
we consider the integral transform
\begin{align}
  \widetilde{M}\left[\lambda^lf((a_i,b_i)_{i=1}^s;\lambda)\right](z)
   &= \frac{1}{2\pi \imu}
       \int_{C}
        \lambda^lf((a_i,b_i)_{i=1}^s;\lambda)z^\lambda\ d\lambda,
    \label{eq.inversMellin}
\end{align}
where $0<z<1$
and the contour $C$ for any $\alpha>0$ is as follows:

 \thinlines
 \unitlength=0.1mm
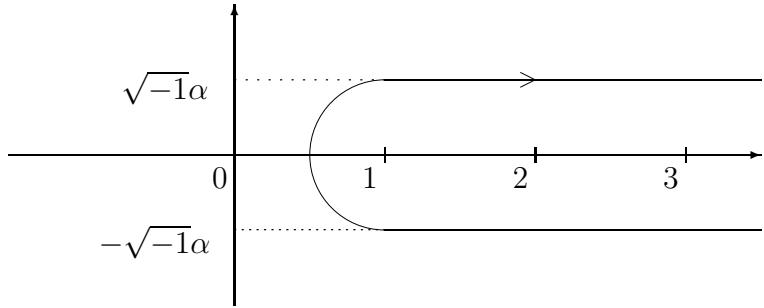
 \begin{picture}(1200,420)(-400,-210)
  \normalfont
  \put(-300,0){\vector(1,0){1000}} 
  \multiput(200,-10)(200,0){3}{\line(0,1){20}}
  \multiputlist(-30,-30)(200,0)[l]{0,1,2,3}
  \put(0,-200){\vector(0,1){400}} 
  \put(700,100){\line(-1,0){500}}
  \put(200,0){\arc{200}{1.5708}{4.7123}}
  \put(200,-100){\line(1,0){500}}
  \put(400,100){\line(-2,1){20}}
  \put(400,100){\line(-2,-1){20}}
  \dottedline{15}(0,100)(200,100)
  \put(-150,75){$\imu\alpha$}
  \dottedline{10}(0,-100)(200,-100)
  \put(-180,-130){$-\imu\alpha$}
 \end{picture}.

\begin{prop}\label{prop.inversMellin}
 The integral transform \eqref{eq.inversMellin} absolutely converges
 and
 \begin{align}
  \widetilde{M}\left[\lambda^l f((a_i,b_i)_{i=1}^s;\lambda)\right](z)
  &=
    \sum_{p=1}^\infty \Res_{\lambda=p}
     \lambda^l f((a_i,b_i)_{i=1}^s;\lambda)\ z^p.
  \label{eq.resphi}
 \end{align}
\end{prop}

\begin{proof}
Since
\begin{align*}
 \left|m_j - \left(t\imu\alpha\right)\right|^2
 =
 (t^2+\alpha^2)\left(1-\frac{m_jt}{t^2+\alpha^2}\right)^2
  + \frac{m_j^2\alpha^2}{t^2+\alpha^2}
 \ge
 \frac{m_j^2\alpha^2}{t^2+\alpha^2},
\end{align*}
we have
\begin{align}
  &\int_0^\infty\left|
   (\pm\imu\alpha + t)^l
    f((a_i,b_i)_{i=1}^s; \pm\imu\alpha + t)
     z^{\pm\imu\alpha +t}\right|dt
   \notag\\
  &\qquad\qquad\le
   \int_0^\infty
    (t^2+\alpha^2)^\frac{l}{2}
    \left|
     f((a_i,b_i)_{i=1}^s; \pm\imu\alpha +t)
    \right|z^{t}\ dt
   \notag\\
  &\qquad\qquad\le
  \zeta(k_1,k_2,\dots,k_n)
   \int_0^\infty z^{t}
  \frac{(t^2+\alpha^2)^{(B_s+l)/2}}{\alpha^{B_s}}
  dt,
  \label{eq.eval2}
\end{align}
where $(k_1,\dots,k_n)=\kappa((a_i,b_i)_{i=1}^s)$
and $B_s=b_1+\dotsb+b_s$.
Thus the integral absolutely converges.

 Next, consider the integral
  \begin{align}
  \frac{1}{2\pi \imu}
   \int_{C_N+\gamma_N}
    \lambda^l
    f((a_i,b_i)_{i=1}^s;\lambda)z^\lambda\ d\lambda
  =
  -\sum_{p=1}^N
   \Res_{\lambda=p} \lambda^l f((a_i,b_i)_{i=1}^s;\lambda)\ z^p
 \end{align}
 where the contour is\\
 \thinlines
 \unitlength=0.1mm
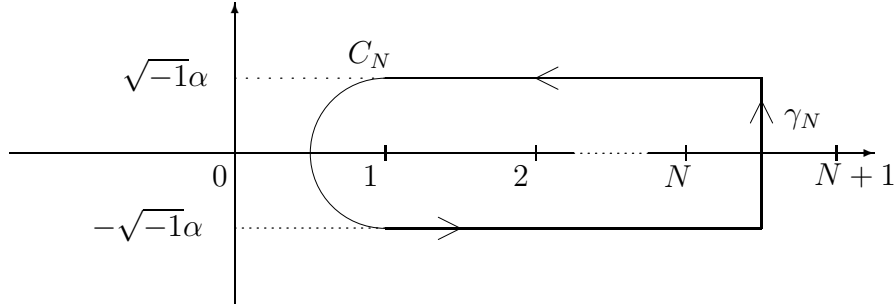
 \begin{picture}(1200,450)(-400,-200)
  \put(-300,0){\line(1,0){750}} 
  \dottedline{10}(450,0)(550,0)
  \put(550,0){\vector(1,0){300}}
   \multiput(200,-10)(200,0){4}{\line(0,1){20}}
   \multiputlist(-30,-30)(200,0)[l]{$0$,$1$,$2$,$N$,$N+1$}
  \put(0,-200){\vector(0,1){400}} 
  \put(200,100){\line(1,0){500}}
   \put(400,100){\line(2,1){30}}
   \put(400,100){\line(2,-1){30}}
  \put(200,0){\arc{200}{1.5708}{4.7123}}
  \put(200,-100){\line(1,0){500}}
   \put(300,-100){\line(-2,1){30}}
   \put(300,-100){\line(-2,-1){30}}
  \put(700,-100){\line(0,1){200}}
   \put(700,70){\line(1,-2){15}}
   \put(700,70){\line(-1,-2){15}}
   \put(730,40){$\gamma_N$}
  \put(150,120){$C_N$}
  \dottedline{15}(0,100)(200,100)
  \put(-150,90){$\imu\alpha$}
  \dottedline{10}(0,-100)(200,-100)
  \put(-190,-110){$-\imu\alpha$}
 \end{picture}\\
 and $\gamma_N$ passes through $N+\frac12$.
 We must prove that
 the integral on $\gamma_N$ tends to $0$ as $N\longrightarrow\infty$.
 Because of the inequality
 \begin{align*}
  \left|\frac{1}{m_j-\left(N+\frac{1}{2}+\imu t\right)}\right|
  \le
  \frac{1}{\left|m_j-\left(N+\frac{1}{2}\right)\right|}
  \le
  \frac{2(N+1)}{m_j},
 \end{align*}
 the integral on $\gamma_N$ is evaluated as
 \begin{align}
  \begin{split}
  &\left|
   \int_{-\alpha}^{\alpha}
    \left(N+\frac{1}{2}+\imu t\right)^l
    f((a_i,b_i)_{i=1}^s;N+\frac{1}{2}+\imu t) z^{N+\frac{1}{2}+\imu t}
   \ dt
  \right|\\
  &\qquad\le
   z^{N+\frac{1}{2}} (2(N+1))^{B_s}
    \zeta(k_1,\dots,k_n)
     \int_{-\alpha}^{\alpha}
      \left(t^2+\left(N+\frac{1}{2}\right)^2\right)^{l/2}\ dt.
  \end{split}
  \label{eq.eval1}
 \end{align}
 Since $0<z<1$, the right hand side converges to $0$.
\end{proof}

We set
\begin{align}
 \varphi((a_i,b_i)_{i=1}^s;z)
   &:= \widetilde{M}\left[f((a_i,b_i)_{i=1}^s;\lambda)\right](z),\\
 \psi((a_i,b_i)_{i=1}^s;z)
   &:= \widetilde{M}\left[g((a_i,b_i)_{i=1}^s;\lambda)\right](z).
\end{align}

\begin{prop}\label{prop.Mellin}
\begin{enumerate}[\upshape (i)]
	\item The functions $\varphi((a_i,b_i)_{i=1}^s;z)$
	      and $\psi((a_i,b_i)_{i=1}^s;z)$ are
	      holomorphic for $|z|<1$.
	\item We have
 \begin{align}
  \widetilde{M}[\lambda^m
  f((a_i,b_i)_{i=1}^s;\lambda)](z)
  &=
  \vartheta^m
  \bigl(\varphi((a_i,b_i)_{i=1}^s;z)\bigr),\\
  \widetilde{M}[(\lambda-1)^mf((a_i,b_i)_{i=1}^s;\lambda-1)](z)
  &=
  z\ \vartheta^m\bigl(\varphi((a_i,b_i)_{i=1}^s;z)\bigr),
 \end{align}
 and
 \begin{align}
 \begin{split}
  &\widetilde{M}
  	\left[
  		\frac{1}{\lambda-1}f((a_i,b_i)_{i=1}^s;\lambda-1)
  	\right](z)\\
  &\qquad\qquad=
 z\left(
  -\zeta(k_1,\dots,k_n)
  + \int_{0}^z \frac{dz}{z}
     \varphi((a_i,b_i)_{i=1}^s;z)
 \right),
 \end{split}
 \end{align}
 where $\vartheta=z d/dz$ is the Euler derivation.
\end{enumerate}
\end{prop}

\begin{proof}
 \begin{enumerate}[$(i)$]
 \item
 From Proposition~\ref{prop.inversMellin} we obtain
 \begin{align}
  \varphi((a_i,b_i)_{i=1}^s;z)
  &=
      \sum_{p=1}^\infty \Res_{\lambda=p}
     f((a_i,b_i)_{i=1}^s;\lambda)\ z^p
    \notag\\
  &=
    \sum_{p=1}^\infty
   \left\{
    \sum_{j=1}^{B_s}
    \sum_{\substack{m_1 > \dots > m_{j-1} > p\\ p >m_{j+1}>\dots>m_{B_s}}}
     C_p^{m_1\dotsc \stackrel{\stackrel{j}\smallsmile}p \dotsc m_{B_s}}\right\}
     z^p,
  \label{eq.seriesphi}
 \end{align}
 where $C_p^{m_1\dotsc \stackrel{\stackrel{j}\smallsmile}p\dotsc m_{B_s}}$
 is the same as in Proposition~\ref{prop.partial}.
 The series
 \begin{align*}
   \sum_{p=1}^\infty
   \left\{
    \sum_{j=1}^{B_s}
    \sum_{\substack{m_1 > \dots > m_{j-1} > p\\ p >m_{j+1}>\dots>m_{B_s}}}
     C_p^{m_1\dotsc \stackrel{\stackrel{j}\smallsmile}p \dotsc m_{B_s}}
   \right\}
   \frac{z^p}{p}
 \end{align*}
 is convergent at $z=1$
 because of Proposition~\ref{prop.partial},
 so the radius of convergence of \eqref{eq.seriesphi} is at least $1$.
 \item
  The first equation can be shown
  by exchanging the derivation and the integration.
  For the second equation
  we shift the integral variable $\lambda$ to $\lambda+1$
  in the left hand side.
 For the last equation, we have
 \begin{align*}
  &\widetilde{M}\left[\frac{1}{\lambda-1}f((a_i,b_i)_{i=1}^s;\lambda-1)\right](z)\\
  &=
  \frac{z}{2\pi\imu}\int_{\{\lambda-1|\lambda\in C \}}
   \frac{1}{\lambda}f((a_i,b_i)_{i=1}^s;\lambda)\ z^\lambda d\lambda\\
  &=
  z
  \left(
   -\zeta(k_1,k_2,\dots,k_n)
   +\frac{1}{2\pi\imu}\int_C
     \frac{1}{\lambda}f((a_i,b_i)_{i=1}^s;\lambda)\ z^\lambda d\lambda
  \right)\\
  &=
  z\left(
   -\zeta(k_1,k_2,\dots,k_n)
   +\int_0^z \left\{
     \frac{1}{2\pi\imu}\int_C
     f((a_i,b_i)_{i=1}^s;\lambda)\ z^\lambda d\lambda
   \right\} \frac{dz}{z}
  \right).
  \end{align*}
  In the last line above we have exchanged
  the order of the integrals.
 \end{enumerate}
\end{proof}

Let us introduce the ``inverse transform'' of $\widetilde{M}$ by
\begin{align}
  M[\varphi(z)](\lambda) = \int_0^1 \varphi(z) z^{-\lambda-1}\ dz.
\end{align}

\begin{prop}
 \begin{align}
  M[\varphi((a_i,b_i)_{i=1}^s;z)](\lambda) 
  =
  f((a_i,b_i)_{i=1}^s;\lambda).
 \end{align}
\end{prop}
\begin{proof}
 For $0<r<1$ and $\lambda < 0$
 \begin{align*}
  &\left|
  \int_0^r 
  \left\{
   \sum_{p=1}^\infty
    \Res_{\lambda=p} f((a_i,b_i)_{i=1}^s;\lambda) z^p
  \right\}
  z^{-\lambda-1}\ dz
  - f((a_i,b_i)_{i=1}^s;\lambda)
  \right|\\
  &\quad=
  \left|
   \sum_{p=1}^\infty
    \Res_{\lambda=p} f((a_i,b_i)_{i=1}^s;\lambda)
     \frac{r^{p-\lambda}}{p-\lambda}
  -
   \sum_{p=1}^\infty
    \Res_{\lambda=p} f((a_i,b_i)_{i=1}^s;\lambda)
     \frac{1}{p-\lambda}
  \right|\\
  &\quad\le
   \sum_{p=1}^\infty
    \left|\Res_{\lambda=p} f((a_i,b_i)_{i=1}^s;\lambda) 
     \frac{1}{p-\lambda}
   \right| (1-r^{p-\lambda})\\
  &\quad\longrightarrow
   0 \qquad (r\to 1)
 \end{align*}
 by virtue of Abel's Theorem.
\end{proof}

\subsection{The differential-integral relations satisfied by
$\varphi$'s and $\psi$'s}

\begin{prop}\label{prop.diffeq}
 The functions $\varphi$'s as well as $\psi$'s
 satisfy the following relations$:$
  \begin{multline}\label{eq.diffeq}
  \sum
   (-\vartheta)^{s-|\delta|-|\epsilon|}
    \varphi((a_i-\delta_i,b_i-\epsilon_i)_{i=1}^s; z)\\
  =
  z \summ 
   (-\vartheta)^{s-|\delta|-|\epsilon|}
    \varphi((a_i-\delta_i,b_i-\epsilon_{i+1})_{i=1}^s; z).
 \end{multline}
 Here $\sum$, $\summ$, $|\delta^{(\prime)}|$, $|\epsilon^{(\prime)}|$
 is the same for Proposition~\ref{prop.sabun}
 and $\vartheta^{-1}$ is the integral operator
 \begin{align}
  \left\{
  \begin{aligned}
   \vartheta^{-1}\varphi\left((a_i,b_i)_{i=1}^s ;z \right)
  &= -\zeta(k_1,\dots,k_n)
   + \int_0^z \frac{dz}{z}
    \varphi
   \left((a_i,b_i)_{i=1}^s;z\right),\\
  \vartheta^{-1}\psi\left((a_i,b_i)_{i=1}^s ;z \right)
  &= -\zeta(k_1',\dots,k_{n'}')
   + \int_0^z \frac{dz}{z}
    \psi
   \left((a_i,b_i)_{i=1}^s;z\right),
  \end{aligned} 
  \right.
  \label{eq.vartheta^{-1}}
\end{align}
 where $(k_1,\dots,k_n)=\kappa((a_i,b_i)_{i=1}^s)$
 and $(k_1',\dots,k_{n'}')$ is its dual.
 We define the weight of $\vartheta$ to be $-1$
 and this relation is homogeneous of weight $|(a_i,b_i)_{i=1}^s| - s$.
\end{prop}
\begin{proof}
 This is a direct consequence from Proposition~\ref{prop.sabun}
 and Proposition~\ref{prop.Mellin}.
\end{proof}

Furthermore we set
(see Definition~\ref{defn.FG} and Proposition~\ref{prop.Mellin})
\begin{align}
 \begin{split}
  \Phi(\bfk;z)
  &:= \widetilde{M}[F(\bfk;z)](z)\\
  &=
   \sum_{\delta_i=0,1} (-\vartheta)^{n-1-|\delta|}
   \varphi((k_1,1)\cup(k_i-\delta_i,1)_{i=2}^n;z),
 \end{split}\\
 \begin{split}
  \Psi(\bfk;z)
  &:= \widetilde{M}[G(\bfk;z)](z)\\
  &=
   \sum_{\delta_i=0,1} (-\vartheta)^{n-1-|\delta|}
   \psi((k_1,1)\cup(k_i-\delta_i,1)_{i=2}^n;z).
 \end{split}
\end{align}
From \eqref{eq.seriesf} and \eqref{eq.resphi} it follows that
\begin{align}
 \Phi(\bfk;z)
 =
 \sum_{m_1>\dots>m_n} \frac{z^{m_1}}{m_1^{k_1}\dotsb m_n^{k_n}}
 =
 \Li{\bfk}{z}.
\end{align}

\begin{cor}\label{cor.Phirel}
The differential relations satisfied by $\Phi(\bfk;z)$ and 
$\Psi(\bfk;z)$ are the same as the differential relations for
$\Li{\bfk}{z}$ 
\begin{align}
 \frac{d}{dz}\Li{k_1,\dots,k_n}{z}
 =
  \begin{cases}
   \dfrac{1}{z}\Li{k_1-1,k_2,\dots,k_n}{z} & \text{if $k_1\ge2$,}\\
   \dfrac{1}{1-z}\Li{k_2,\dots,k_n}{z}   & \text{if $k_1=1$}.
  \end{cases}
\end{align}
\end{cor}

\begin{proof}
This is clear from Proposition~\ref{prop.Frel} and 
Proposition~\ref{prop.Mellin}.
\end{proof}

\subsection{Mellin transform and inverse Mellin transform}

We call the integral transforms $M$ and $\widetilde{M}$
Mellin transform and inverse Mellin transform respectively.
In fact, for suitable compositions,
$\widetilde{M}$ is
actually inverse Mellin transform.

\begin{prop}
For any composition $(a_i,b_i)_{i=1}^s$
with $a_{i_0} \ge 2$ for some $i_0$ and $0<c<1$
\begin{align}
  \widetilde{M}\left[f((a_i,b_i)_{i=1}^s;\lambda)\right](z)
   = \frac{1}{2\pi \imu}\int_{c-\imu\infty}^{c+\imu\infty}
      f((a_i,b_i)_{i=1}^s;\lambda)z^\lambda\ d\lambda.
\end{align}
\end{prop}

\begin{proof}
 Consider the following contour

 \thinlines
 \unitlength=0.08mm
 \begin{picture}(1200,620)(-400,-310)
  \put(-300,0){\line(1,0){750}} 
  \dottedline{10}(450,0)(550,0)
  \put(550,0){\vector(1,0){300}}
   \multiput(200,-10)(200,0){4}{\line(0,1){20}}
   \multiputlist(-30,-30)(200,0)[l]{$0$,$1$,$2$,$N$,$N+1$}
  \dottedline{15}(700,-100)(700,100)
  \put(0,-300){\vector(0,1){600}} 
  \put(200,100){\line(1,0){500}}
   \put(400,100){\line(2,1){30}}
   \put(400,100){\line(2,-1){30}}
  \put(200,0){\arc{200}{1.5708}{4.7123}}
  \put(200,-100){\line(1,0){500}}
   \put(300,-100){\line(-2,1){30}}
   \put(300,-100){\line(-2,-1){30}}
  \put(150,120){$C_N$}
  \put(700,-200){\line(0,1){100}}
   \put(700,-170){\line(1,2){15}}
   \put(700,-170){\line(-1,2){15}}
   \put(730,-170){$\gamma_1^-$}
  \put(0,-200){\line(1,0){700}}
   \put(300,-200){\line(2,1){30}}
   \put(300,-200){\line(2,-1){30}}
   \put(350,-240){$\gamma_2^-$}
   \put(-140,-220){$\scriptstyle -\imu L$}
  \put(0,0){} 
   \put(0,50){\line(1,-2){15}}
   \put(0,50){\line(-1,-2){15}}
  \put(0,200){\line(1,0){700}}
   \put(-110,200){$\scriptstyle \imu L$}
   \put(200,200){\line(-2,1){30}}
   \put(200,200){\line(-2,-1){30}}
   \put(350,220){$\gamma_2^+$}
  \put(700,100){\line(0,1){100}}
   \put(700,140){\line(-1,2){15}}
   \put(700,140){\line(1,2){15}}
   \put(730,140){$\gamma_1^+$}
  \dottedline{15}(0,100)(200,100)
  \put(-100,90){$\scriptstyle\imu\alpha$}
  \dottedline{10}(0,-100)(200,-100)
  \put(-130,-110){$\scriptstyle-\imu\alpha$}
 \end{picture}.

 Because the integrand has singular points only at positive integers
 the integral on this contour is $0$.
 We check that the integrals
 on $\gamma_1^+$, $\gamma_1^-$, $\gamma_2^+$
 and $\gamma_2^- \to 0$
 if $N$, $L \to \infty$.
 First let N to $\infty$,
 then we can show that the integral on $\gamma_1^+$ tends to 0
 in the same manner as in \eqref{eq.eval1}.
 Similarly the integral on $\gamma_1^- \longrightarrow$  0.

 Next let L to $\infty$. Assume that $i_0=1$.
 Noting that
 \begin{align*}
  \frac{1}{\left| m_1^{a_1} -\left(t+\imu L\right) \right|}
  \le \frac{1}{m_1^{a_1}L}
 \end{align*}
 we can show that in the same manner as in \eqref{eq.eval2}
 the integral on $\gamma_2^+$ is
 \begin{align*}
  &\left|\int_0^\infty z^{\imu L+t} f((a_i,b_i)_{i=1}^s;\imu L+t)\ dt\right|\\
  &\qquad\qquad\le
  \frac{\zeta(k_1-1,k_2,\dots,k_n)}{L}
   \int_0^\infty z^{t}
  \left(\frac{t^2+L^2}{L^2}\right)^{B_s-1}
  \ dt\\
  &\qquad\qquad \longrightarrow 0\qquad(L\longrightarrow\infty).
 \end{align*}
 The other cases of $i_0 \neq 1$ and that the integral 
 on $\gamma_2^-$ tends to $0$ can be verified in the same way.
\end{proof}

\section{Proof and examples of Theorem~\ref{thm.redohno}}
\label{sec.proofredohno}
\subsection{Proof of Theorem~\ref{thm.redohno}}

From the equation \eqref{eq.diffeq}
in Proposition~\ref{prop.diffeq}
we obtain
\begin{multline*}
 \left\{
  (1-z)\sum_{\substack{\epsilon_j\delta_{j+1}=0\\\forall j}}
  +\sum_{
   \substack{\epsilon_j\delta_{j+1}=1\\{}^\exists j}}
 \right\}
  (-\vartheta)^{s-|\delta|-|\epsilon|}
   \varphi((a_i-\delta_i,b_i-\epsilon_i)_{i=1}^s; z)\\
 =
 z\summ_{\substack{\delta_j \epsilon_{j+1}=1\\ \exists j}}
  (-\vartheta)^{s-|\delta|-|\epsilon|}
   \varphi((a_i-\delta_i,b_i-\epsilon_{i+1})_{i=1}^s; z).
\end{multline*}
Here we note that
$\varphi((a_i,b_i)_{i=1}^s; z)$ appears only in the first summation
(it corresponds to the term with all $\delta_i=\epsilon_i=0$).
Divding by $\frac{1}{1-z}$
and applying $\int_0^z\frac{dz}{z}$ $s$ times
we have
\begin{multline*}
  (-)^s\varphi((a_i,b_i)_{i=1}^s; z)
  -
  \sum_{\substack{
 		\epsilon_j\delta_{j+1}=0\\\forall j\\
		1\le|\delta|+|\epsilon|\le s}
       }
  \left(-\int_0^z\frac{dz}{z}\right)^{s-|\delta|-|\epsilon|}
   \varphi((a_i-\delta_i,b_i-\epsilon_i)_{i=1}^s; z)\\
  +
  \sum_{
   \substack{\epsilon_j\delta_{j+1}=1\\{}^\exists j}}
  \left(\int_0^z\frac{dz}{z}\right)^{s}
  (-\vartheta)^{s-|\delta|-|\epsilon|}
   \varphi((a_i-\delta_i,b_i-\epsilon_i)_{i=1}^s; z)\\
  +
  \sum_{
   \substack{\epsilon_j\delta_{j+1}=1\\{}^\exists j}}
  \left(\int_0^z\frac{dz}{z}\right)^{s-1}
  \int_0^z \frac{dz}{1-z}
  (-\vartheta)^{s-|\delta|-|\epsilon|}
   \varphi((a_i-\delta_i,b_i-\epsilon_i)_{i=1}^s; z)\\
 =
  \summ_{\substack{\delta_j \epsilon_{j+1}=1\\ \exists j}}
  \left(\int_0^z\frac{dz}{z}\right)^{s-1}
 \int_0^z
 \frac{dz}{1-z}
  (-\vartheta)^{s-|\delta|-|\epsilon|}
   \varphi((a_i-\delta_i,b_i-\epsilon_{i+1})_{i=1}^s; z).
\end{multline*}
From the induction hypothesis
$\varphi$'s of less compositions than $(a_i,b_i)_{i=1}^s$
can be written by MZVs-linear combination of $\Phi$'s.
Using Corollary~\ref{cor.Phirel}
the second and third terms in the left hand side
are expressed as MZVs-linear combination of $\Phi$'s.
 So we have
 \begin{multline*}
  \varphi((a_i,b_i)_{i=1}^s)
  = \sum_\bfc \alpha_\bfc^{(a_i,b_i)} \zeta(\bfk_\bfc^{(a_i,b_i)})
        \Phi(\bfc;z)\\
  +
   \sum_\bfc \beta_\bfc^{(a_i,b_i)} \zeta(\bfk_\bfc^{(a_i,b_i)})
    \left(\int_0^z\frac{dz}{z}\right)^{s-1}
    \int_0^z\frac{dz}{1-z} \vartheta^{s-m_\bfc}
    \bigl(\Phi(\bfc;z)\bigr),
 \end{multline*}
 where $2\le m_\bfc \le s+1$
 and $\alpha_\bfc^{(a_i,b_i)}$, $\beta_\bfc^{(a_i,b_i)}\in\Q$.
 For the proof we must verify that
 the last term above is written in MZV-linear combination of $\Phi$'s.
 If $m_\bfc=s$, from Corollary~\ref{cor.Phirel}
 \begin{align*}
     \left(\int_0^z\frac{dz}{z}\right)^{s-1}
    \int_0^z\frac{dz}{1-z} \vartheta^{s-m_\bfc}
    \bigl(\Phi(\bfc;z)\bigr)
   =
     \left(\int_0^z\frac{dz}{z}\right)^{s-1}
    \int_0^z\frac{dz}{1-z}
    \bigl(\Phi(\bfc;z)\bigr)
 \end{align*}
 is $\Phi$ of the greater composition.
 If $m_\bfc=s+1$, one can show by virtue of \eqref{eq.vartheta^{-1}},
 \begin{align*}
  \left(\int_0^z\frac{dz}{z}\right)^{s-1}
   \int_0^z\frac{dz}{1-z} \vartheta^{-1}
    \bigl(\Phi(\bfc;z)\bigr).
 \end{align*}
 is written in MZV-linear combination of $\Phi's$.
 For the $\Psi$'s case the corresponding coefficients are
 MZVs of the dual indices.
 If $2 \le m_\bfc \le s-1$,
 $\vartheta^{s-m_\bfc}\bigl(\Phi(\bfc;z)\bigr)$
 can be expressed as linear combination
 of the following 
 \begin{align}
  \vartheta^{s-m_\bfc}
    \bigl(\Phi(\bfc;z)\bigr)
  =
  \sum_{\bfc',l,m}
   \alpha^{(a_i,b_i)}
   \Phi(\bfc',z) \frac{z^l}{(1-z)^m},
  \label{eq.dels}
 \end{align}
 where $\bfc'$ is less than $\bfc$,
 $l\le m \le s-m_\bfc$
 and $\alpha^{(a_i,b_i)}\in\Z$.
 Next, we repeat the iterated integrals
 using the expansions
\begin{xalignat*}{2}
  \frac{z^l}{(1-z)^m}
  &= \sum_{j=0}^{n}\binom{n}{j}\frac{(-1)^{n-j}}{(1-z)^{m-n+j}},&
  \frac{1}{(1-z)^mz}
  &= \frac{1}{z} + \sum_{j=1}^m \frac{1}{(1-z)^j}.
\end{xalignat*}
Making once the iterated integral of
the right hand side of \eqref{eq.dels}
we have the following:
\begin{align*}
 \begin{split}
  &\int_0^z\frac{dz}{(1-z)^mz}\Phi(\bfk;z)\\
  &\qquad=
   \int_0^z\frac{dz}{z} \Phi(\bfk;z)
   +\int_0^z\frac{dz}{1-z} \Phi(\bfk;z)
   +\sum_{j=2}^m \int_0^z \frac{dz}{(1-z)^j}\Phi(\bfk;z),
 \end{split}
\end{align*}
\begin{align*}
 &\int_0^z \frac{dz}{(1-z)^m}\Phi(k_1,\dots,k_n;z)\\
 &=
 \frac{1}{m-1}
 \left\{
  \frac{1}{(1-z)^{m-1}}\Phi(k_1,\dots,k_n;z)
   - \Phi(k_1,\dots,k_n;z)
 \right.\\
  &\left.
   \qquad
   - \Phi(1,k_1-1,\dots,k_n;z)
   - \sum_{l=2}^{m-1}\int_0^z \frac{dz}{(1-z)^{l}}\Phi(k_1-1,\dots,k_n;z)
  \right\},
\end{align*}
\begin{align*}
 \begin{split}
   &\int_0^z \frac{dz}{(1-z)^m}\Phi(1,k_2,\dots,k_n;z)\\
  &=
  \frac{1}{m-1}
  \left\{
   \frac{1}{(1-z)^{m-1}}\Phi(1,k_2,\dots,k_n;z)
    - \int_0^z \frac{dz}{(1-z)^m}\Phi(k_2,\dots,k_n;z)
  \right\},
  \end{split}
\end{align*}
and
\begin{align*}
 \int_0^z \frac{dz}{(1-z)^m}
 = -z\sum_{l=1}^{m-1}\frac{1}{l}\binom{m-1}{l-1}\frac{z^{l-1}}{(1-z)^{m-1}}.
\end{align*}
Hence after making the iterated integral $s-1$ times
we finally reach to the integral
\begin{align*}
 \int_0^z \frac{dz}{1-z} = \Phi(1;z).
\end{align*}

 Thus 
 $\varphi((a_i,b_i)_{i=1}^s)$ can be written
 by MZVs-linear combination
 of $\Phi$'s.
 The application of Mellin transform
 \begin{align*}
  M[\varphi(z)](\lambda) =
  \int_0^1 \varphi(z)z^{-\lambda-1} dz
 \end{align*}
 gives us the theorem.
\qed

\subsection{Examples of Theorem~\ref{thm.redohno}}

We list examples of \eqref{eq.fasF} up to weight $6$.
For simplicity we drop the variable $\lambda$.

weight $2$:
\begin{align}
 f(1,1) &= F(1).
\end{align}

weight $3$:
\begin{align}
 f(2,1) &= f(1,2) = F(2).
\end{align}

weight $4$:
\begin{align}
 f(3,1) = f(1,3) &= F(3),\\
 f(2,2) &= 2 F(3) + F(1,2) - \zeta(2)F(1),\\
 f(1,1,1,1) &= F(3)-F(2,1).
\end{align}

weight $5$:
\begin{align}
 f(4,1) &= f(1,4) = F(4),\\
 f(3,2) &= 3F(4) + F(1,3) + F(2,2)
   - \zeta(3)F(1) - \zeta(2)F(2),\\
 f(2,3) &= 3F(4) + F(1,3) + F(2,2)
   -\zeta(2,1)F(1) - \zeta(2)F(2),\\
 f(2,1,1,1) &= f(1,1,1,2)
   = 2F(4) - F(3,1) + F(2,2) - \zeta(2) F(2),\\
 f(1,2,1,1) &= f(1,1,2,1)
  = F(4) - F(3,1) - F(2,2).
\end{align}

weight $6$:
\begin{align}
 &f(5,1) = f(1,5) = F(5),\\
 \begin{split}
  &f(4,2) = 4F(5) + F(3,2) + F(2,3) + F(1,4)\\
  & \qquad\qquad\qquad
    - \zeta(4)F(1) - \zeta(3)F(2) - \zeta(2)F(3),
 \end{split}\\
 \begin{split}
  &f(2,4) = 4F(5) + F(3,2) + F(2,3) + F(1,4)\\
   &\qquad\qquad\qquad
              - \zeta(2,1,1)F(1) - \zeta(2,1)F(2)
              - \zeta(2)F(3),
 \end{split}\\
 \begin{split}
  &f(3,3) =
  6F(5) + 2F(3,2) + 2 F(2,3) + 2 F(1,4) + F(1,2,2)\\
  &\qquad\qquad
  - \zeta(3,1)F(1) 
  - \zeta(3)F(2) - \zeta(2,1)F(2)
  - 2\zeta(2)F(3) - \zeta(2)F(1,2),
 \end{split}\\
 &f(3,1,1,1) = 3F(5) -F(4,1) + F(2,3) + F(3,2)
   - \zeta(3)F(2) - \zeta(2)F(3),\\
 &f(1,1,1,3) = 3F(5) -F(4,1) + F(2,3) + F(3,2)
   - \zeta(2,1)F(2) - \zeta(2)F(3),\\
 \begin{split}
  &f(2,2,1,1) = 3F(5) - 2F(4,1) - F(2,2,1) - F(2,1,2)\\
  &\qquad\qquad\qquad\qquad\qquad\qquad
   -\zeta(2)F(3) + \zeta(2)F(2,1) - \zeta(2,1)F(2),
 \end{split}\\
 \begin{split}
  &f(1,1,2,2) = 3F(5) - 2F(4,1) - F(2,2,1) - F(2,1,2)\\
  &\qquad\qquad\qquad\qquad\qquad\qquad
   -\zeta(2)F(3) + \zeta(2)F(2,1) - \zeta(3)F(2),
 \end{split}\\
 &f(2,1,2,1) = f(1,2,1,2)
  = 2F(5) - F(4,1) - \zeta(2)F(3),\\
 \begin{split}
  &f(2,1,1,2) = 5F(5) - F(4,1) + 3F(3,2) + 2F(2,3) + F(2,1,2)\\
  &\qquad\qquad\qquad\qquad\qquad\qquad\qquad\qquad\qquad
  - 3\zeta(2)F(3) - \zeta(2)F(2,1),
 \end{split}\\
 &f(1,3,1,1) = f(1,1,3,1) = F(5) - F(4,1) - F(3,2) - F(2,3),\\
 &f(1,2,2,1) = F(5) - F(4,1) - 2F(3,2) -F(2,3) + F(2,2,1),\\
 &f(1,1,1,1,1,1) = 2F(5) -2F(4,1) + F(3,1,1) - \zeta(2)F(3).
\end{align}

\section{Landen formula and the Ohno relation}
\label{sec.landen}

The Landen connection formula for the dilogarithm \eqref{eq.landilog}
generalizes to the MPLs case.
This formula is interpreted as the connection formula
of MPLs between $1$ and $\infty$.
\begin{prop}\label{prop.landen}
 \begin{align}
 \Li{\bfk}{z} = (-1)^n \sum_{\bfc_1,\dots,\bfc_n}
 \Li{\bfc_1\dotsm\bfc_n}{\frac{z}{z-1}},
  \label{eq.landen again}
 \end{align}
 where $\bfc_i$ runs all compositions of $k_i$.
\end{prop}

\begin{proof}
 For $\bfk=(1)$
 \begin{align*}
  \Li{1}{z}= -\log(1-z) 
   = \log\left(1-\frac{z}{z-1}\right)
   = -\Li{1}{\frac{z}{z-1}}.
 \end{align*}

 We assume that the proposition holds
 for $\bfk=(k_1,k_2,\dots,k_n)$.
 Then using the differential relation
 \begin{align}
  \begin{split}
   \frac{d}{dz}\Li{k_1,\dots,k_n}{\frac{z}{z-1}}
   =
   \begin{cases}
    \left(\dfrac{1}{z}+\dfrac{1}{1-z}\right)
    \Li{k_1-1,k_2,\dots,k_n}{\dfrac{z}{z-1}}
    & \text{if $k_1\ge2$,}\\
    -\dfrac{1}{1-z}\Li{k_2,\dots,k_n}{\dfrac{z}{z-1}}
    & \text{if $k_1=1$}
   \end{cases}
  \end{split}\label{eq.diffeq4landen}
 \end{align}
 we have
 \begin{align*}
 & \Li{k_1+1,k_2,\dots,k_n}{z}
 =
  \int_0^z \frac{dz}{z} \Li{k_1,k_2,\dots,k_n}{z}\\
 &=
 \int_0^z
 \left(
  \left(\frac{dz}{z}+\frac{dz}{1-z}\right)dz - \frac{dz}{1-z}
  \right)
  (-1)^n \sum_{\bfc_1,\dots,\bfc_n}
   \Li{\bfc_1\dotsm \bfc_n}{\frac{z}{z-1}},
 \end{align*}
 and
 \begin{align*}
  \Li{1,k_2,\dots,k_n}{z}
 &=
  \int_0^z \frac{dz}{1-z} \Li{k_1,k_2,\dots,k_n}{z}\\
 &=
 -\int_0^z -\frac{dz}{1-z}
  (-1)^n \sum_{\bfc_1,\dots,\bfc_n}
   \Li{\bfc_1\dotsm\bfc_n}{\frac{z}{z-1}}.
 \end{align*}
\end{proof}

One can obtain further information of
the right hand side of \eqref{eq.landen again}
in the case of $\bfk=k$, a positive integer.

\begin{lem}\label{lem.landen}
 For any positive integers $j$ and $k$ with $j\le k$,
 we have
 \begin{align}
  \sum_{m_1>\dots>m_k}
  \frac{z^{m_j}}{m_1\prod_{i\ne j}(m_i-m_j)}
  =
  -\sum_{\substack{\bfc\\ |\bfc|=k\\ \len(\bfc) = k-j+1}}
   \Li{\bfc}{\frac{z}{z-1}}.
 \end{align}
\end{lem}

\begin{proof}
 We show by induction for $k$.
 For the case $k=1$
 \begin{align*}
  \sum_{m=1}^\infty \frac{z^m}{m} = - \Li{1}{\frac{z}{z-1}}
 \end{align*}
 is obvious.
 We suppose that the proposition is correct for $k-1$.
 Calculating the derivative of the series
 and applying the induction hypothesis
 we have
 \begin{align*}
 \begin{split}
 &\frac{d}{dz}
  \sum_{m_1>\dots>m_k}
  \frac{z^{m_j}}{m_1\prod_{i\ne j}(m_i-m_j)}\\
  &=
 \frac{1}{z(1-z)}
  \sum_{m_1>\dots>m_{k-1}}
  \frac{z^{m_{j-1}}}{m_1\prod_{i\ne j}(m_i-m_{j-1})}\\
 &\hspace{5cm}
 -
 \frac{1}{1-z}
  \sum_{m_1>\dots>m_{k-1}}
  \frac{z^{m_j}}{m_1\prod_{i\ne j}(m_i-m_j)}
 \end{split}\\
  &=
  -\frac{1}{z(1-z)}
  \sum_{\substack{\bfc\\ |\bfc|=k-1\\ \len(\bfc) = k-j+1}}
  \Li{\bfc}{\frac{z}{z-1}}
  +\frac{1}{z}
  \sum_{\substack{\bfc\\ |\bfc| = k-1 \\ \len(\bfc) = k-j}}
  \Li{\bfc}{\frac{z}{z-1}}.
 \end{align*}
 Due to the differential relation \eqref{eq.diffeq4landen}
 the lemma is proved for $k$.
\end{proof}

For any positive integer $k$, we have
\begin{align}
 G(k;\lambda) &=
 \sum_{m_1>\dots>m_k}
  \frac{1}{m_1(m_1-\lambda)\dotsb(m_k-\lambda)},
 \label{eq.G_k}
  \\
\intertext{hence}
 \Psi(k;z) &=
   \sum_{m_1>\dots>m_k}
 \sum_{j=1}^k
  \frac{z^{m_j}}{m_1\prod_{i\ne j}(m_i-m_j)}.
 \label{eq.Psi_k}
\end{align}
By virtue of Lemma~\ref{lem.landen} we have
\begin{prop}
 \begin{align}
 \Psi(k;z) =
  -\sum_{j=1}^k
  \sum_{\substack{\bfc\\ |\bfc| = k \\ \len(\bfc) = k-j+1}}
   \Li{\bfc}{\frac{z}{z-1}}.
 \end{align}
 Moreover using the differential relation for Li's and $\Psi$'s
 and the equation above,
 we can see that
 \begin{align}
 \Psi(k_1,\dots,k_n;z)
 = (-1)^n \sum_{\bfc_1,\dots,\bfc_n}
 \Li{\bfc_1\dotsm \bfc_n}{\frac{z}{z-1}}.
 \end{align}
\end{prop}
Thus the relation $\Phi=\Psi$ can be interpreted
as the Landen connection formula.
We can think that
the reduced Ohno relation is converted,
via inverse Mellin transform,
to the connection formula of MPLs between $1$ and $\infty$.

\noindent\textbf{Discussions}.
\quad Proposition~\ref{prop.landen} says that the
reduced Ohno relation is converted to the Landen connection
formula for MPLs by inverse Mellin transform.
In particular, since the explicit forms \eqref{eq.G_k}, \eqref{eq.Psi_k}
of $G(k;\lambda)$ and $\Psi(k;z)$ are revealed,
the Landen connection formula for the polylogarithm $\Li{k}{z}$
for any positive integer $k$ is, via Mellin transform, converted
to the relation $F(k;\lambda)=G(k;\lambda)$. This is nothing
but the sum formula of depth $k$. Thus the sum formulas for MZVs
are equivalent to the Landen connection formulas for
polylogarithms. However, such equivalency is not achieved yet
for the reduced Ohno relation of indices $\bfk$ of
depth greater than $1$.
This is an important issue to be settled in future.

We have another important issue: In Theorem~\ref{thm.redohno},
which is the main theorem in this paper, the indices
$\bfk_{\bfc}^{(a_i,b_i)}$
and the rational numbers $\alpha_{\bfc}^{(a_i,b_i)}$
are not specified. They must be determined. Moreover, we
conjecture that
\begin{align}
    f((a_i,b_i)_{i=1}^s;\lambda) =&
     \sum_{\bfc}\alpha_{\bfc}^{(a_i,b_i)}
            \zeta(\bfk_{\bfc}^{(a_i,b_i)})
                \ast   F(\bfc;\lambda),  \\
      g((a_i,b_i)_{i=1}^s;\lambda) =&
     \sum_{\bfc}\alpha_{\bfc}^{(a_i,b_i)}
            \zeta(\bfk_{\bfc}^{\prime (a_i,b_i)})
                \ast   G(\bfc;\lambda),
\end{align}
where $\ast$ is the harmonic product introduced in \cite{hoff2}.

\section{Proofs of the difference relations and the Ohno relation}
\label{sec.proofsabun}

\subsection{Proof of Proposition~\ref{prop.sabun}}

To prove Proposition~\ref{prop.sabun}, we need the following.

\begin{lem}\label{lem.sabun}
 For any composition $(a_i,b_i)_{i=1}^s$, we set
 \begin{multline}
   [\{(a_i,d_i),b_i\}_{i=1}^s;\lambda]\\
   := \sum_{m_1 > \dots > m_{B_s} > 0}
   \prod_{i=1}^s
   \frac{1}{(m_{B_{i-1}+1}-d_i)^{a_i}
    \underbrace{
     (m_{B_{i-1}+1}-\lambda)\dotsm(m_{B_i}-\lambda)
    }_{b_i}
   },
 \end{multline}
where we interpret special cases with $a_i=0$ or $b_i=0$ for some $i$
as follows{\normalfont :}
 \begin{align}
   [\{\ldots,b_{i-1},(0,d_i),b_i,\ldots\};\lambda]
   &= [\{\ldots,b_{i-1}+b_i,\ldots\};\lambda], \\
   [\{\ldots,(a_{i-1},d),0,(a_i,d),\ldots\};\lambda]
   &= [\{\ldots,(a_{i-1}+a_i,d),\ldots\};\lambda].
 \end{align}
Then we have the following difference relations$:$
\begin{enumerate}[{\normalfont \upshape (i)}]
 \item 
 \begin{enumerate}[{\normalfont \upshape (a)}]
  \item\label{lemma.2} If $a_1\ge2$,
   \begin{align}
    \begin{split}
      &\lambda\,[\{(a_1,0),b_1,\dots\};\lambda]\\
     &\qquad\quad
     - [\{(a_1-1,0),b_1,\dots\};\lambda]
     - [\{(a_1,0),b_1-1,\dots\};\lambda]\\
     &=
     \lambda'\,[\{(a_1,1),b_1,\dots\};\lambda]
      - [\{(a_1-1,1),b_1,\dots\};\lambda].
    \end{split}   
   \end{align}
  \item\label{lem.sabun.ib} If $a_1=1$,
   \begin{multline}
     \lambda\,[\{(1,0),b_1,\dotsc\};\lambda]
      - [\{(1,0),b_1-1,\dotsc\};\lambda]\\
    = \lambda'\,[\{(1,1),b_1,\dotsc\};\lambda].
   \end{multline}
 \end{enumerate}
 \item\label{lem.sabun.ii}
       If $i\not=1$ and $i\not=s$, or $i=s$ and $b_s\not=1$,
  \begin{align}
   \begin{split}
    &\lambda\,[\{\dots,(a_i,0),b_i,\dots\};\lambda] \\
    &\qquad
    - [\{\dots,(a_i-1,0),b_i,\dots\};\lambda]
    - [\{\dots,(a_i,0),b_i-1,\dots\};\lambda]\\
    &=
     \lambda'\,[\{\dots,b_{i-1},(a_i,1),\dots\};\lambda]\\
     &\quad
      - [\{\dots,b_{i-1},(a_i-1,1),\dots\};\lambda]
   - [\{\dots,b_{i-1}-1,(a_i,1),\dots\};\lambda].
   \end{split}
  \end{align}
 \item
  \begin{enumerate}[{\normalfont \upshape (a)}]
   \item\label{lemma.5} If $b_s \ge 2$,
  \begin{multline}
   [\{(a_1,1),b_1,(a_2,1),b_2,\dots,(a_s,1),b_s\};\lambda]\\
    = [\{(a_1,0),b_1,\dots,(a_s,0),b_s\};\lambda']\\
    - \frac{1}{\lambda'} [\{(a_1,0),b_1,\dots,(a_s,0),b_s-1\};\lambda'].
  \end{multline}
 \item If $b_s=1$,
  \begin{align}
   \begin{split}
    &\lambda\,[\{(a_1,1),b_1,(a_2,1),b_2,\dots,
      (a_{s-1},1),b_{s-1},(a_s,0),1\};\lambda]\\
    &- [\{(a_1,1),b_1,(a_2,1),b_2,\dots,
      (a_{s-1},1),b_{s-1},(a_s-1,0),1\};\lambda]\\
    &=
    \lambda'\,[\{(a_1,0),b_1,(a_2,0),b_2,\dots,(a_s,0),1\};\lambda']\\
    & \qquad
    - [\{(a_1,0),b_1,(a_2,0),b_2,\dots,(a_s-1),1\};\lambda']\\
    & \qquad\quad
     -[\{(a_1,0),b_1,(a_2,0),b_2,\dots,b_{s-1}-1,(a_s,0),1\};\lambda'].
   \end{split}  
  \end{align} 
 \end{enumerate}
\end{enumerate}
\end{lem}

\begin{proof}
We use the partial-fractions expansion:
\begin{align}
 \begin{split}
  &\frac{\lambda}{m^a(m-\lambda)} - \frac{1}{m^{a-1}(m-\lambda)}\\
  &\quad= 
  \frac{\lambda'}{(m-1)^a(m-\lambda)} - \frac{1}{(m-1)^{a-1}(m-\lambda)}
  + \left(\frac{1}{(m-1)^a} - \frac{1}{m^a} \right),
 \end{split}\label{frac.1}\\
 &\frac{\lambda}{m(m-\lambda)}
 = \frac{\lambda'}{(m-1)(m-\lambda)}
 + \left(\frac{1}{m-1} - \frac{1}{m}\right).
 \label{frac.2}
\end{align}

\begin{enumerate}[{\normalfont \upshape (i)}]
\item
\begin{enumerate}[{\normalfont \upshape (a)}]
\item
We set $B$ by
\begin{align*}
 \prod_{j=1}^s
   m_{B_{j-1}+1}^{a_j}
    \underbrace{
     (m_{B_{j-1}+1}-\lambda)\dotsm(m_{X_j}-\lambda)
    }_{b_j}
   = m_1^{a_1}(m_1-\lambda)\,B.
\end{align*}
Then using \eqref{frac.1} we have
 \begin{align*}
  &\lambda\,[\{(a_1,0),b_1,\dotsc\};\lambda]
    - [\{(a_1-1,0),b_1,\dotsc\};\lambda]\\
  &\qquad=
   \sum_{m_1 > \dots>m_{B_s} > 0}
   \left\{
    \frac{\lambda}{m_1^{a_1}(m_1-\lambda)}
    - \frac{1}{m_1^{a_1-1}(m_1-\lambda)}
   \right\}\frac{1}{B}\\
 \begin{split}
   &\qquad=
   \sum_{m_1 > \dots>m_{B_s} > 0}
   \Biggl\{
    \frac{\lambda'}{(m_1-1)^{a_1}(m_1-\lambda)}
    - \frac{1}{(m_1-1)^{a_1-1}(m_1-\lambda)}\\
    & \hspace{5,5cm}
    + \left(\frac{1}{(m_1-1)^{a_1}}-\frac{1}{m_1^{a_1}}\right)
   \Biggr\}\frac{1}{B}
 \end{split} \\
 \begin{split}
   &\qquad=
   \sum_{m_1 > \dots>m_{B_s} > 0}
   \left\{
    \frac{\lambda'}{(m_1-1)^{a_1}(m_1-\lambda)}
    - \frac{1}{(m_1-1)^{a_1-1}(m_1-\lambda)}
   \right\}\frac{1}{B}\\
   &\hspace{4cm}
    +\sum_{m_2>\dots>m_{B_s}>0}\,\sum_{m_1 = m_2+1}^\infty
     \left(\frac{1}{(m_1-1)^{a_1}}-\frac{1}{m_1^{a_1}}\right)\frac{1}{B}
 \end{split}\\
  &\qquad=
  \begin{aligned}[t]
   &\lambda'\,[\{(a_1,1),b_1,\dotsc\};\lambda]
    -[\{(a_1-1,1),b_1,\dotsc\};\lambda]\\
   &\hspace{4cm}
    +[\{(a_1,0),b_1-1,\dotsc\};\lambda].
  \end{aligned}
 \end{align*}

\item
 Using \eqref{frac.2}, it can be proved in the same manner as \eqref{lemma.2}.
\end{enumerate}
\item We set $A$ and $B$ by
  \begin{align*}
   &A := \prod_{j=1}^{i-1}
   m_{B_{j-1}+1}^{a_j}
    \overbrace{
     (m_{B_{j-1}+1}-\lambda)\dotsm(m_{B_j}-\lambda)
    }^{b_j}, \\
   &B:=
   \frac
    {\prod_{j=i-1}^s
    m_{B_{j-1}+1}^{a_j}
     \overbrace{
      (m_{B_{j-1}+1}-\lambda)\dotsm(m_{B_{j}}-\lambda)
     }^{b_j}}
    {m_{B_{i-1}+1}^{a_i}(m_{B_{i-1}+1}-\lambda)}.
  \end{align*}
Then
 \begin{align*}
  &\lambda [\{\dots,b_{i-1},(a_i,0),b_i,\dots\};\lambda]
   - [\{\dots,b_{i-1},(a_i-1,0),b_i,\dots\};\lambda]\\
  &\quad=
   \begin{aligned}[t]
    &\sum_{m_1 > \dots > m_{B_s} > 0}
    \frac{1}{A}
    \left\{
    \frac{\lambda}{
                    m_{B_{i-1}+1}^{a_i}(m_{B_{i-1}+1}-\lambda)}
    -\frac{1}{
                 m_{B_{i-1}+1}^{a_i-1}(m_{B_{i-1}+1}-\lambda)}
    \right\} \frac{1}{B}\\
   \end{aligned}\\
  &\quad=
   \sum_{m_1 > \dots > m_{B_s} > 0}
   \frac{1}{A}
    \Biggl\{
  \begin{aligned}[t]
   &  \frac{\lambda'}{(m_{B_{i-1}+1}-1)^{a_i}(m_{B_{i-1}+1}-\lambda)}\\
   &\qquad
      -\frac{1}{(m_{B_{i-1}+1}-1)^{a_i-1}(m_{B_{i-1}+1}-\lambda)}\\
   &\qquad\qquad
   \left.-\left(
     \frac{1}{(m_{B_{i-1}+1}-1)^{a_i}} - \frac{1}{m_{B_{i-1}+1}^{a_i}}
   \right) 
   \right\}\frac{1}{B}
  \end{aligned}\\
  &\quad=
  \begin{aligned}[t]
   &\lambda' [\dots,b_{i-1},(a_i,1),b_i,\dots;\lambda]
    - [\dots,b_{i-1},(a_i-1,1),b_i,\dots;\lambda]\\
   &\quad
   + 
   \!\!\!
   \sum_{\substack{m_1> \dots > m_{B_{i-1}}\\
            m_{B_{i-1}} -1 >m_{B_{i-1}+2} > \dots > m_{B_s} > 0}}
   \!\!\!
   \frac{1}{A}
   \left(
    \frac{1}{(m_{B_{i-1}+2})^{a_i}}
     - \frac{1}{(m_{B_{i-1}}-1)^{a_i}}
   \right) \frac{1}{B}.
  \end{aligned}
 \end{align*}
We divide the range of sum of the third term into two parts as
  \begin{align*}
   \sum_{\substack{m_1> \dots > m_{B_{i-1}}\\
    m_{M_{i-1}} -1 >m_{B_{i-1}+2} > \dots > m_{B_s} > 0}}
   = 
   \sum_{n_1> \dots > m_{B_s} > 0}
    - \sum_{\substack{m_1> \dots > m_{B_{i-1}}\\
       m_{B_{i-1}+2} = m_{B_{i-1}} -1 \\
       m_{B_{i-1}+2} > \dots > m_{B_s} > 0}}.
  \end{align*}
The later sum is equal to zero because of $m_{B_{i-1}+2} = m_{B_{i-1}} -1$.
Thus we have
 \begin{align*}
  &\lambda [\{\dots,b_{i-1},(a_i,0),b_i,\dots\};\lambda]
   - [\{\dots,b_{i-1},(a_i-1,0),b_i,\dots\};\lambda]\\
  &=
  \begin{aligned}[t]
   &\lambda' [\{\dots,b_{i-1},(a_i,1),b_i,\dots\};\lambda]
    - [\{\dots,b_{i-1},(a_i-1,1),b_i,\dots\};\lambda]\\
   & + \sum_{m_1> \dots > m_{B_s} > 0}
   \frac{1}{A}
   \left\{
    \frac{1}{(m_{B_{i-1}+2})^{a_i}}
     - \frac{1}{(m_{B_{i-1}}-1)^{a_i}}
   \right\} \frac{1}{B}
  \end{aligned}\\
  &=
  \begin{aligned}[t]
    \lambda'& [\{\dots,b_{i-1},(a_i,1),b_i,\dots\};\lambda]
    - [\{\dots,b_{i-1},(a_i-1,1),b_i,\dots\};\lambda]\\
    &- [\{\dots,b_{i-1}-1,(a_i,1),b_i,\dots\};\lambda]
    +[\{\dots,b_{i-1},(a_i,0),b_i-1,\dots\};\lambda].
  \end{aligned}
 \end{align*}
\item
\begin{enumerate}[{\normalfont \upshape (a)}]
\item[(b)]
 Repeating shift of $m_i \mapsto m_i+1$, we have
 \begin{align*}
  &\begin{aligned}[t]
   &\lambda\,[\{(a_1,1),b_1,
      \dots,(a_{s-1},1),b_{s-1},(a_s,0),1\};\lambda]\\
   &\qquad
   -\lambda\,[\{(a_1,1),b_1,
      \dots,(a_{s-1},1),b_{s-1},(a_s-1,0),1\};\lambda]
  \end{aligned}\\
  &\qquad=
  -\sum_{m_1>\dots>m_{B_s}>0}
   \frac{1}{(m_1-1)^{a_1}\dotsm (m_{B_s-1}-\lambda)m_{B_s}^{a_s}}\\
  \begin{split}
   &\qquad=
  -\sum_{m_1>\dots>m_{B_s}\ge 0}
   \frac{1}{m_1^{a_1}\dotsm (m_{B_s-1}-\lambda')(m_{B_s}+1)^{a_s}}\\
  &\qquad\qquad\text{(by shift $m_i\longmapsto m_i+1$)}
  \end{split}\\
  \begin{split}
   &\qquad=
  -\sum_{m_1>\dots>m_{B_s-1}\ge m_{B_s}> 0}
   \frac{1}{m_1^{a_1}\dotsm (m_{B_s-1}-\lambda')m_{B_s}^{a_s}}\\
  &\qquad\qquad\text{(by shift $m_{B_s}+1 \longmapsto m_{B_s}$)}
  \end{split}\\
  \begin{split}
   &\qquad=
  -\sum_{m_1>\dots>m_{B_s}> 0}
   \frac{m_{B_s}-\lambda'}
    {m_1^{a_1}\dotsm(m_{B_s-1}-\lambda')
      m_{B_s}^{a_s}(m_{B_s}-\lambda')}\\
  &\qquad\qquad\qquad\qquad
   -\sum_{m_1>\dots>m_{B_s-1}> 0}
   \frac{1}{m_1^{a_1}\dotsm m_{B_s-1}^{a_s}(m_{B_s-1}-\lambda')}
  \end{split}\\
   \begin{split}
    &\qquad=
    \lambda'\,[\{(a_1,0),b_1,\dots,b_{s-1}-1,(a_s,0),1\};\lambda'] \\
    &\qquad\qquad
     - [\{(a_1,0),b_1,\dots,b_{s-1},(a_s-1,0),1\};\lambda']\\
    &\qquad\qquad\qquad\qquad
     - [\{(a_1,0),b_1,\dots,b_{s-1}-1,(a_s,0),1\};\lambda'].
   \end{split} 
 \end{align*} 
\item
Similarly as in the previous cases,
 \begin{align*}
  &[\{(a_i,1),b_i\}_{i=1}^s;\lambda]\\
  &\quad=
   \sum_{m_1>\dots>m_{B_s}>0}
    \frac{1}{(m_1-1)^{a_1}(m_1-\lambda)
             \dots (m_{B_s-1}-\lambda)(m_{B_s}-\lambda)}\\
  &\quad=
   \sum_{m_1>\dots>m_{B_s}\ge 0}
    \frac{1}{m_1^{a_1}(m_1-\lambda')
             \dots (m_{B_s-1}-\lambda')(m_{B_s}-\lambda')}\\
  &\quad=
   (\text{$m_{B_s}>0$ part}) + (\text{$b_{B_s}=0$ part})\\
  &\quad=
   [\{(a_i,0),b_i\}_{i=1}^s;\lambda']
    - \frac{1}{\lambda'}[\{(a_i,0),b_i\}_{i=1}^{s-1}\cup\{(a_s,0),b_s-1\};\lambda'].
 \end{align*}
\end{enumerate}
\end{enumerate}
\end{proof}

\begin{proof}[\bf Proof of Proposition~\ref{prop.sabun}]
Using $[\{(a_i,d_i),b_i\}_{i=1}^s]$,
the generating functions $f$ and $g$ are expressed as follows:
\begin{align*}
 \left\{
 \begin{aligned}[c]
  &f((a_i,b_i)_{i=1}^s;\lambda) = [\{(a_i,0),b_i\}_{i=1}^s;\lambda], \\
  &g((a_i,b_i)_{i=1}^s;\lambda) = [\{(b_i,0),a_i\}_{i=s}^1;\lambda].
 \end{aligned}
 \right.
\end{align*} 
If $a_1,b_s\ge2$,
applying Lemma~\ref{lem.sabun} successively
 \begin{align*}
 &\text{(LHS)}\\
 &=
  \sum_{\{(\delta_i,\epsilon_i)\}_{i=1}^s\in I^s}
   (-\lambda)^{s-|\delta|-|\epsilon|}\,
   [\{(a_i-\delta_i,0),b_i-\epsilon_i\}_{i=1}^s;\lambda]\\
\begin{split}
 &=
   \sum_{\{(\delta_i,\epsilon_i)\}_{i=2}^s\in I^{s-1}}
  \sum_{\delta_1'\in\{0,1\}}
    (-\lambda)^{s-1-|\delta|-|\epsilon|}(-\lambda')^{1-|\delta'|}\\
    &\qquad\times
     [\{(a_1-\delta_1',1),b_1,(a_2-\delta_2,0),b_2-\epsilon_2,
      \dots,
      (a_s-\delta_s,0),b_s-\epsilon_s\};\lambda]\\
  &\qquad\qquad\text{(by Lemma~\ref{lem.sabun}~(\ref{lemma.2}))}
\end{split} \\
  \begin{split}
   &=
   \sum_{\delta_1'\in\{0,1\}}
    \sum_{\{(\delta_j',\epsilon_j')\}_{j=2}^s\in I^{s-1}}
     (-\lambda')^{s-|\delta'|-|\epsilon'|}\\
      &\quad\times
       [\{(a_1-\delta_1',1),b_1-\epsilon_2',
        (a_2-\delta_2',1),b_2-\delta_3',
         \dots,
          (a_s-\delta_s',1),b_s\};\lambda] \\
          &\qquad\qquad
   \text{(by Lemma~\ref{lem.sabun}~(\ref{lem.sabun.ii}) $s-1$ times)}
  \end{split}\\
  \begin{split}
   &=
    \sum_{\delta_1',\epsilon_{s+1}'\in\{0,1\}}
     \sum_{\{(\delta_j',\epsilon_j')\}_{j=2}^s\in I^{s-1}}
      (-\lambda')^{s-|\delta'|-|\epsilon'|}\\
       &\qquad\times
        [\{(a_1-\delta_1',0),b_1-\epsilon_2',(a_2-\delta_2',0),b_2-\delta_3',\\
         &\hspace{6cm}\dots,
          (a_s-\delta_s',0),b_s-\epsilon_{s+1}'\};\lambda']\\
    &\qquad\qquad\text{(by Lemma~\ref{lem.sabun}~(\ref{lemma.5}))}
  \end{split}\\
&=
  \text{(RHS).}
 \end{align*} 
Remaining relations
  and the relations of $g$'s can be proved quite similarly.
\end{proof}

\subsection{Alternative proof of the Ohno relation}

 From the properties of the generating functions
 clarified in Section~\ref{sec.diffrel},
 we give an alternative proof for the Ohno relation
 \begin{align*}
  f((a_i,b_i)_{i=1}^s;\lambda) =  g((a_i,b_i)_{i=1}^s;\lambda)
 \end{align*}
 by induction on compositions.

 If the composition is minimum i.e. $(a_i,b_i)_{i=1}^s = (1,1)$,
 it is obvious.

 If the theorem is correct for compositions
 less than $(a_i,b_i)_{i=1}^s$,
 applying Proposition~\ref{prop.sabun} to
 $(a_i,b_i)_{i=1}^s$ for $f$ and $g$,
 we obtain two relations for $f$'s and $g$'s.
 Subtracting these two equations, we have
 \begin{align*}
  &\sum{\lambda}^{s-|\delta| - |\epsilon|}\,
   \Bigl\{
   f((a_i-\delta_i,b_i-\epsilon_i)_{i=1}^s;\lambda)
   - g((a_i-\delta_i,b_i-\epsilon_i)_{i=1}^s;\lambda)
  \Bigr\}
  \\
  &\quad=
  \sum{\lambda'}^{s-|\delta'|-|\epsilon'|}
   \Bigl\{
    f((a_i-\delta_i',b_i-\epsilon_i')_{i=1}^s;\lambda')
   - g((a_i-\delta_i',b_i-\epsilon_i')_{i=1}^s;\lambda')
   \Bigr\}.
 \end{align*}
 But the terms whose compositions are less than $(a_i,b_i)_{i=1}^m$
 are canceled out
 by the induction hypothesis.
 The remaining is
 \begin{multline*}
  {\lambda}^s
   \Bigl\{
    f((a_i,b_i)_{i=1}^s;\lambda)
   - g((a_i,b_i)_{i=1}^s;\lambda)
   \Bigr\}
   \\
  \qquad=
  {\lambda'}^s
  \Bigl\{
   f((a_i,b_i)_{i=1}^s;\lambda')
   - g((a_i,b_i)_{i=1}^s;\lambda')
  \Bigr\}.
 \end{multline*}
 Hence
  ${\lambda}^s\,f((a_i,b_i)_{i=1}^s;\lambda)
    - {\lambda}^s\,g((a_i,b_i)_{i=1}^s;\lambda) $
 is a periodic function in $\lambda$ with period $1$.
 Furthermore by Proposition~\ref{prop.partial} it is a meromorphic
 function such as
 \begin{align*}
  \lambda^s \sum_{p=1}^\infty
   \frac{C_k}{p-\lambda}.
 \end{align*}
 Because of the periodicity, all $C_k$'s must be zero.
 Thus we complete the proof.

\end{document}